\documentclass[12pt]{amsart}
\usepackage{amsmath,amssymb}

\begin{document}
\newtheorem{thm}{Theorem}[section]
\newtheorem{lem}[thm]{Lemma}

\newtheorem{cor}[thm]{Corollary}
\newtheorem{conj}[thm]{Conjecture}
\newtheorem{clm}[thm]{Claim}
\theoremstyle{definition}
\newtheorem{dfn}[thm]{Definition}
\newtheorem{notation}[thm]{Notation}
\newtheorem{question}[thm]{Question}
\theoremstyle{remark}

\newtheorem{exm}[thm]{Example}
\newtheorem{rem}[thm]{Remark}
\def\N{{\mathbb N}}
\def\G{{\mathbb G}}
\def\Q{{\mathbb Q}}
\def\R{{\mathbb R}}
\def\C{{\mathbb C}}
\def\P{{\mathbb P}}
\def\A{{\mathbb A}}
\def\Z{{\mathbb Z}}
\def\v{{\mathbf v}}
\def\x{{\mathbf x}}
\def\O{{\mathcal O}}
\def\M{{\mathcal M}}
\def\kbar{{\bar{k}}}
\def\tr{\mbox{Tr}}
\def\id{\mbox{id}}
\def\qed{{\tiny $\clubsuit$ \normalsize}}

\renewcommand{\theenumi}{\alph{enumi}}

\title{The arithmetic puncturing problem and integral points}

\author{David McKinnon}
\author{Yi Zhu}

\begin{abstract}
In this paper, we provide evidence to support a positive answer to a question of Hassett and Tschinkel.  In particular, if an algebraic variety $V$ has a dense set of rational points, they ask whether or not the set of $D$-integral points is potentially dense, where $D$ is a set of codimension at least two.  We give a positive answer to this question in many cases, including varieties whose generic linear section is a smooth rational curve, and certain $K3$ surfaces.  We also discuss some stronger notions of integrality of points, and give some positive answers to some cases of the analogous question in the stronger context.
\end{abstract}

\maketitle

\vspace{-5pt}
\section{Introduction}

In \cite{HT}, as Problem 2.13 (``The Arithmetic Puncturing Problem"), Hassett and Tschinkel ask the following question:  

\begin{question}\label{htquestion}
Let $X$ be a projective variety with canonical singularities and $D$ a subvariety of codimension $\geq 2$.  Assume that rational points on $X$ are potentially dense.  Are integral points on $(X,D)$ potentially dense?
\end{question}

\vspace{.1in}

In the paper, they provide positive answers to this question in various cases, including toric varieties and products of elliptic curves.  The purpose of this paper is to provide a larger set of positive answers to Question~\ref{htquestion}, to ask a stronger version of Question~\ref{htquestion} (namely Question~\ref{mzquestion}), and to give a positive answer to the stronger question in a large number of cases.

Of course, the hypothesis that $D$ has codimension two cannot be removed, as there are countless well known examples of varieties with a dense set of rational points but a degenerate set of integral points if $D$ is a divisor.  

In order to state the stronger question, we need to fix some notation and definitions.  Let $X$ be a projective algebraic variety, $D\subset X$ a Zariski closed subset, both defined over a number field $k$.  Let $M_k$ be the set of places of $k$.  The following definition is a slight generalization of Definition 1.4.3 in \cite{Vo}:

\begin{dfn}
Let $S$ be a finite set of places of $k$.  A subset $R\subset X(k)-D(k)$ is called {\em $(D,S)$-integralizable} if and only if there are global Weil functions $\lambda_{D,v}$ and non-negative real numbers $n_v$ such that $n_v=0$ for all but finitely many places $v$ for each $v$, and such that
\[\lambda_{D,v}(P)\leq n_v\]
for all $v\in M_k-S$ and $P\in R$.
\end{dfn}

This generalizes \cite{Vo} only in that we do not demand that $S$ contain the set of infinite places of $k$.  For definitions and discussion of Weil functions (called ``local height functions'' by \cite{Si}), see \cite{La}, \cite{Si}, and \cite{Vo}.  

The particular case of this new definition that most interests us in this paper is the following:

\begin{dfn}
A subset $R\subset X(k)$ is {\em everywhere $D$-integral} if and only if $R$ is $(D,\emptyset)$-integralizable.
\end{dfn}

We can now state the stronger version of Question~\ref{htquestion}:

\begin{question}\label{mzquestion}
Let $X$ be a projective variety with canonical singularities and $D$ a subvariety of codimension $\geq 2$.  Assume that rational points on $X$ are potentially dense.  Are everywhere $D$-integral points on $X$ potentially dense?
\end{question}

The classical notion of an integral point is simple one whose coordinates don't have denominators.  In affine space over $\Q$, an integral point in this sense is one whose closure in a model over $\mbox{Spec}(\Z)$ does not intersect the closure of the hyperplane at infinity.  One can generalize this notion to that of an $S$-integral point, by permitting intersections over a finite set $S$ of places.  If $S'$ is the union of $S$ with the set of infinite places of $k$, then this corresponds precisely to our notion of $(D,S')$-integralizable, only with the additional requirement that $n_v=0$ for all finite $v\not\in S$.

The remainder of the paper is structured as follows.  Section 2 establishes a few preliminaries, regarding various different notions of integral point and comparing them.  Section 3 contains the main theorems regarding everywhere integral points and Question~\ref{mzquestion}, and Section 4 proves results about more classical notions of integral points and Question~\ref{htquestion}.  The central philosophy of the paper is to prove preliminary results on curves, and then use those to prove results in higher dimension when the variety in question is covered by curves.  

\section{Funding}

The first author was supported by a Discovery Grant from the Natural Sciences and Engineering Research Council of Canada.

\section{Preliminaries}

In this paper, we will consider four different notions of integrality:

\begin{enumerate}
\item Total $(D,\emptyset)$-integrality, where $n_v=0$ for all finite places $v$ but there are still constraints at the infinite places.  
\item  Sets of integral points in the truly classical sense, which correspond to $(D,\infty)$-integral points, where $\infty$ represents the set of infinite places of $k$.
\item Everywhere $D$-integral sets.
\item $(D,S)$-integralizable sets.
\end{enumerate}  

Note that the last two notions of integrality apply only to sets -- they are not well defined for individual points.  Applied to sets, each of these four notions of integrality are related in the following ways:

Classical $(D,\emptyset)$-integrality is the strongest notion of integrality that we will consider, but Zariski dense sets of classical $(D,\emptyset)$-integral points do exist.  For example, if $D$ is empty, then every rational point is classically $(D,\emptyset)$-integral.  Any set of classically $(D,\emptyset)$-integral points is also $(D,\infty)$-integral, everywhere $D$-integral, and $(D,S)$-integralizable for any $S$.

The notions of ``classically $(D,\infty)$-integral sets of points'' and ``everywhere $D$-integral sets'' are incomparable.  The set of points in $\P^1(\Q)$ of the form $[n:1]$ for $n\in\Z$ is classically $D$-integral for $D=\infty$, but is not everywhere $D$-integral, because if $\lambda_\infty$ is the Weil function for $\infty$ at the infinite place, then for any real number $B$, there are only finitely many $n$ for which $\lambda_\infty(n:1)\leq n$.  Conversely, let $X=\P^2$, with $D=\{[0:0:1],[0:1:0],[0:1:1],[1:0:0],[1:0:1],[1:1:0],[1:1:1]\}$.  Then since every rational point in $\P^1$ is congruent to one of the seven points in $D$ modulo $2$, there are no classically $D$-integral points ... but by setting $n_2=2$, one easily obtains an infinite set of everywhere $D$-integral points, for example points of the form $\{[2k+1:2:2]\}$ for any odd integer $k$.

Finally, a set which satisfies any of the first three definitions of integrality is automatically $(D,S)$-integral for any choice of $S$ (if we require $S$ to contain all the infinite places).

We would like to be able to prove results for classical everywhere $D$-integral sets of points, but have not yet been able to surmount the additional technical obstacles created by the additional restriction.

If $D$ is a very ample divisor, then there are no infinite, everywhere $D$-integral sets.  

\begin{thm}
Let $X$ be a projective algebraic variety, $D$ a very ample divisor on $X$ (all defined over a number field $k$), $R$ an everywhere $D$-integral subset of $X(k)$.  Then $R$ is finite.
\end{thm}

\noindent
{\it Proof:} \/ We start by proving the theorem for $X=\P^n_k$ and $D=\{x_0=0\}$.  Any set of everywhere $D$-integral points is a set of $k$-rational points in the affine piece $x_0\neq 0$ whose coordinates have $v$-adic absolute value absolutely bounded for each infinite place $v$.  Moreover, there is an integer $N$ such that every point in $R$ can be written as a fraction of integers in $\O_k$ whose denominator divides $N$.  Thus, every coordinate of every point of the set $R$ lies in a compact subset of Minkowski space (that is, the product of the archimedean completions of $k$ modulo complex conjugation), and is a subset of the lattice $\frac{1}{N}\O_k$, so by the geometry of numbers $R$ must be finite.  It follows that the set $R$ must be finite.

More generally, let $X$ be an arbitrary projective variety, $D$ an arbitrary very ample divisor.  Then there is an embedding $\phi\colon X\to\P^n$ such that $\phi^*(\{x_0=0\})=D$.  Furthermore, $\phi(R)$ must be an everywhere $\{x_0=0\}$-integral set, so it must be finite, implying the finiteness of $R$.  \qed

\vspace{.1in}

By contrast, if $D=\emptyset$, then an everywhere $D$-integral set is simply a set of rational points, and there are many examples of Zariski dense sets of rational points.  In this article we provide nontrivial examples of such sets, in order to illuminate the spectrum of possibilities.

\section{Everywhere Integral Points}

We begin with a theorem about everywhere integral points on curves with respect to an arithmetic zero-cycle.

\begin{thm}\label{curveintegral}
Let $V$ be an algebraic variety defined over a number field $k$, and let $\mathcal{V}$ be a model of $V$ over $\mbox{Spec}(\O_k)$.  Let $C$ be an irreducible curve on $V$, and let $\mathcal{C}$ be its closure in $\mathcal{V}$.  Let $D$ be a Zariski closed subset of $V$, and let $\mathcal{D}$ be its closure in $\mathcal{V}$.

For every place $v$ of $k$, let $n_v$ be a non-negative real number such that $n_v=0$ for all but finitely many places $v$.  Choose Weil functions $\lambda_{D,v}$ for each place $v$.  Assume that there is a point $P\in C(k)$ satisfying
\[\lambda_{D,v}(P)\leq n_v\]
for every place $v$.

If $D\cap C=\emptyset$ and $C(k)$ is infinite, then there are infinitely many points $Q\in C(k)$ satisfying
\[\lambda_{D,v}(Q)\leq n_v\]
\end{thm}

\noindent
%
%

%
%
%

\noindent
{\it Proof:} \/ Let $S$ be the set of places $v$ of $k$ such that either $v$ is archimedean, or $\mathcal{D}\cap\mathcal{C}$ is supported on $v$.  Note that $S$ is finite.

For each $v\not\in S$, $\lambda_{D,v}(Q)=0$ for all $Q\in C$, so we may restrict our attention to $v\in S$.

If $v$ is finite with corresponding prime $\pi$ of $\O_k$, then the condition $\lambda_{D,v}(P)<n_v$ depends only on the residue class of $P$ modulo a suitable power of $\pi$.  (See for example subsection 2.2.2 of \cite{BG}.)  Thus, the collection of all $Q\in C(k)$ satisfying $\lambda_{D,v}(Q)\leq n_v$ for all finite $v$ contains the set of points $Q$ such that $Q\equiv P\pmod N$ for some suitable nonzero $N\in\O_k$.

If $v$ is infinite, then the restriction $\lambda_{D,v}(Q)\leq n_v$ defines the complement of an open disk in the $v$-adic topology.

There are now two cases: either the geometric genus of $C$ is zero, or one.

If the geometric genus of $C$ is zero, then the theorem now follows immediately from the Weak Approximation Theorem for $\P^1$.

Weak Approximation does not hold for curves of genus one, however, so we must work a bit harder.  The set $B$ of points $Q$ such that $Q\equiv P\pmod N$ for some nonzero $N\in \O_k$ contains the image on $C$ of a coset of a finite index subgroup $A$ of the Mordell-Weil group of the normalization $\tilde{C}$ of $C$ over $k$.  It remains to show that $B$ contains an infinite subset satisfying $\lambda_{D,v}(Q)\leq n_v$ for all the infinite places $v$.

The subgroup $A$ is infinite by assumption, so it contains a point $R\neq P$ of infinite order.  If we let $\C'$ be the completion of $\tilde{C}$ with respect to the product of all the infinite $v$-adic topologies, it becomes a finite union of complex tori.  By Kronecker's Theorem in diophantine approximation (see for example Theorem 443 of \cite{HW}), the identity element of $C'$ is a limit point of the subgroup $A$.  After translation by some element in the inverse image of $R$, and pushing forward to $C$, we see that $B$ has infinite intersection with arbitrarily small neighbourhoods of $R$ (in the product of infinite topologies).  A small enough such neighbourhood will consist entirely of points satisfying $\lambda_{D,v}(Q)\leq n_v$ for all the infinite places $v$, so the theorem is proven.  \qed

\vspace{.1in}


The main theorems of this section crucially use Theorem~\ref{curveintegral} to establish the density of integral points in higher dimension.

\begin{thm}\label{rcvaintegral}
Let $X$ be a projective variety defined over a number field $k$, embedded in $\P^n$ so that the intersection of $X$ with a generic linear space of codimension $\dim X-1$ is a curve of genus zero.  Let $D$ be an algebraic subset of $X$ of codimension at least two.  Then, possibly after some finite extension of $k$, there is a Zariski dense set of everywhere $D$-integral points on $X$.  In particular, Question~\ref{mzquestion} has a positive answer for $(X,D)$.
\end{thm}

\noindent
{\it Proof:} \/ We first choose a model $\mathcal{X}$ for $X$ over $\mbox{Spec}\O_k$, and let $\mathcal{D}$ be the closure of $D$ in $\mathcal{X}$.  Choose a global Weil function $\lambda_D$ for $D$, with corresponding local Weil functions $\lambda_{D,v}$.

Choose an algebraic point $P$ not contained in the support of $D$ or the singular set of $X$, and define non-negative real numbers $n_v$ by $n_v=\lambda_{D,v}(P)$.  After a finite extension of $k$, we may assume that $P$ is $k$-rational.  Since this is the only time when a finite extension of $k$ is necessary, we will fix that finite extension now and assume hereafter that $k$ is large enough.  

A generic linear space of codimension $\dim X-1$ through $P$ intersects $X$ in a curve of genus zero, which must be birational to $\P^1_k$ by Bertini's Theorem because it contains the smooth point $P$.  Let $A$ be the set of linear spaces $L$ that satisfy the following conditions:
\begin{itemize}
\item $L$ is defined over $k$.
\item $L$ contains $P$.
\item $L\cap X$ is a curve of genus zero that is smooth at $P$.
\item $L\cap D$ is empty, so that $\mathcal{L}\cap\mathcal{D}$ is an arithmetic zero-cycle, where $\mathcal{L}$ is the closure of $L$ in $\mathcal{X}$.
\end{itemize}
By Theorem~\ref{curveintegral}, for any $L\in A$, there are infinitely many rational points $Q$ on $L\cap X$ such that
\[\lambda_{D,v}(Q)\leq n_v\]
for all $v$.  The union of such curves is dense in $X$.  To see this, note that the union of all linear spaces $L$ that pass through $P$ and are defined over $k$ is Zariski dense in $X$.  The conditions of being smooth at $P$ and intersecting $D$ trivially are both nonempty open conditions (because $D$ has codimension at least two in $X$), so the union of the curves in $A$ must be dense in $X$.  We conclude that there is a Zariski dense set of everywhere $D$-integral points, as desired.  \qed

\vspace{.1in}

\begin{rem}
	Theorem~\ref{rcvaintegral} applies in a great many cases.  For example, if $X$ is a cone over a rational normal curve in projective space and $D$ is a finite set of points on $X$, then Theorem~\ref{rcvaintegral} implies that the set of everywhere $D$-integral points on $X$ is potentially Zariski dense.  
\end{rem}

We can prove a variation on Theorem~\ref{rcvaintegral} using the group structure on an elliptic surface.

\begin{thm}\label{ellintegral}
Let $\pi\colon X\to\P^1$ be an elliptic surface, defined over a number field $k$, and let $\mathcal{E}$ be the corresponding elliptic curve over $k(t)$.  Let $S$ be the identity section of $\pi$ (that is, the section corresponding to the identity element of $\mathcal{E}$).  Assume that there is a curve $C$ on $X$ with an infinite set of $k$-rational points, and such that $C$ is not a component of a fibre of $\pi$, such that $C$ does not correspond to a torsion point of $\mathcal{E}$.  Let $D$ be a finite set of points on $X$ such that $D\cap S=\emptyset$.  Then there is an everywhere $D$-integral set of points on $X$.  In particular, Question~\ref{mzquestion} has a positive answer for $(X,D)$.
\end{thm}	

\noindent
{\it Proof:} \/ By replacing $C$ with a suitable multiple of itself (in the group law on $\mathcal{E}$), we may assume that $C\cap D=\emptyset$.  (For each fibre $F$ of $\pi$ that intersects $D$, either $C\cap F$ has infinite order on $F$ -- in which case any large enough multiple of $C$ will be disjoint from $D$ on $F$ -- or else multiply $C$ by a positive integer $n$ divisible enough so that $nC\cap F=S\cap F\not\in D$.)  
Let $\mathcal{X}$ be a regular model of $X$ over $\mbox{Spec}(\O_k)$, and let $\mathcal{C}$ and $\mathcal{D}$ be the closures of $C$ and $D$, respectively, in $\mathcal{X}$.  

By Theorem~\ref{curveintegral}, $C$ contains an infinite, everywhere $D$-integralizable set of points.

Moreover, by Silverman's Specialization Theorem (Theorem 11.4 of Chapter III of \cite{Si2}), we may assume that all but finitely many of those everywhere $D$-integral points lies on a smooth fibre of $\pi$, and has infinite order as an element of that fibre.

We conclude that there is an infinite set $G$ of $k$-rational points on $\P^1$ such that for all $P\in G$, $\pi^{-1}(P)$ is a smooth elliptic curve containing an everywhere $D$-integral point of infinite order, and such that $\pi(D)\cap G=\emptyset$.  

Any fibre $F$ of $\pi$ with $\pi(F)\in G$ is therefore a smooth elliptic curve with $F\cap D=\emptyset$, and with an everywhere $D$-integral point $P$ of infinite order.  Again by Theorem~\ref{curveintegral}, all such fibres contain an infinite set of everywhere $D$-integral points.  \qed

\section{$S$-integral Points}

We begin by proving the following variation on Theorem~\ref{curveintegral}.  It will play a similarly central role in this section.

\begin{thm}\label{curvesintegral}
Let $V$ be an algebraic variety defined over a number field $k$, and let $\mathcal{V}$ be a model of $V$ over $\mbox{Spec}(\O_k)$.  Let $C$ be an irreducible curve on $V$, and let $\mathcal{C}$ be its closure in $\mathcal{V}$.  Let $D$ and $N$ be Zariski closed subsets of $V$, and let $\mathcal{D}$ and $\mathcal{N}$ be their closures in $\mathcal{V}$, respectively.  Let $L=D\cup N$, and let $S$ be a set of places of $k$ that contains all the archimedean places of $k$.

For every place $v$ of $k$ with $v\not\in S$, let $n_v$ be a non-negative real number such that $n_v=0$ for all but finitely many places $v$.  Choose Weil functions $\lambda_{L,v}$ for each place $v$.  Assume that there is a point $P\in C(k)$ satisfying
\[\lambda_{L,v}(P)\leq n_v\]
for every place $v\not\in S$.

If $N\cap C=\emptyset$ and $C$ contains an infinite set of $(D,S)$-integral points, then there are infinitely many points $Q\in C(k)$ satisfying
\[\lambda_{L,v}(Q)\leq n_v\]
for every place $v\not\in S$.
\end{thm}

\noindent
{\it Proof:} \/  Since $C(k)$ is infinite, it follows that $C$ must have geometric genus zero or one.  The condition that $C$ contain a dense set of $D$-integral points implies that $C$ must intersect $D$ in at most two places of $C$ (places in the sense of points of the normalization of $C$), and that $C\cap D=\emptyset$ if $C$ has genus one.

The case in which $C\cap D=\emptyset$ follows {\em a fortiori} from Theorem~\ref{curveintegral}, so we assume that $C$ is a genus zero curve with either one or two places supported on $D$.  Let $\pi\colon \tilde{C}\to C$ be the normalization map over $k$.  The set of points $R$ of $\tilde{C}$ with $\pi(R)\not\in D$ are a principal homogeneous space for an arithmetic group ($\mathbb{G}_a$ if there is one place of $C$ on $D$, and $\mathbb{G}_m$ if there are two places), so by choosing a point $R_0$ on $\tilde{C}$, we can give the $(\pi^*D,S)$-integral points of $\tilde{C}$ the structure of an arithmetic group $G$.  The set $B$ of points $Q$ of $C$ such that $Q\equiv P\pmod N$ for some nonzero $N\in\O_k$ -- which, as before, is contained in the set of points $Q$ satisfying $\lambda_{L,v}(Q)\leq n_v$ for all $v$ not in $S$ -- contains the image on $C$ of a coset of a finite index subgroup $A$ of $G$, and is therefore infinite, as desired.  \qed

\begin{rem}
Note that the naive analogue of Theorem~\ref{curvesintegral} in the setting of Theorem~\ref{curveintegral} is false.  That is, if $C$ is isomorphic to $\mathbb{G}_a$ or $\mathbb{G}_m$, then even if $C\cap D$ is empty (over $k$), it's possible that the set of everywhere $D$-integral points of $C$ are finite and nonempty.  In both cases, the integral points of $C$ are a finitely generated abelian group.  If we choose $k=\Q$, and consider the archimedean place of $\Q$, then the only limit points of the integral points of $C$ are the place(s) at infinity.  Thus, if the set $\{Q\mid \lambda_{D,\infty}(Q)\leq n_\infty\}$ does not include a place at infinity, then it will be finite and discrete.
\end{rem}

Theorem~\ref{curvesintegral} enables us to prove the following result, in the same spirit as Theorem~\ref{rcvaintegral}.

\begin{thm}
Let $X$ be a projective variety defined over a number field $k$, embedded in $\P^n$ so that the intersection of $X$ with a generic linear space of codimension $\dim X-1$ is a curve of genus zero.  Let $D$ be a proper Zariski closed subset of $X$ whose codimension one (in $X$) component has degree at most two.  Let $S$ be a finite set of places of $k$, containing all the archimedean places.  Then, possibly after some finite extension of $k$, there is a Zariski dense set of everywhere $(D,S)$-integral points on $X$.
\end{thm}

\noindent
{\it Proof:} \/ The proof is eerily similar to the proof of Theorem~\ref{rcvaintegral}.  We first choose a model $\mathcal{X}$ for $X$ over $\mbox{Spec}\O_k$, and let $\mathcal{D}$ be the closure of $D$ in $\mathcal{X}$.  Choose a global Weil function $\lambda_D$ for $D$, with corresponding local Weil functions $\lambda_{D,v}$.  Also, write $D=D'+N$, where $D'$ is a (Weil) divisor and $N$ is a cycle of codimension at least two.  (The components of $N$ need not all have the same dimension.)

Choose an algebraic point $P$ not contained in the support of $D$ or the singular set of $X$, and define non-negative real numbers $n_v$ by $n_v=\lambda_{D,v}(P)$.  After a finite extension of $k$, we may assume that $P$ is $k$-rational, and that the unit group of the ring of integers $\O_k$ is infinite.  Since this is the only time when a finite extension of $k$ is necessary, we will fix that finite extension now and assume hereafter that $k$ is large enough.  

A generic linear space of codimension $\dim X-1$ through $P$ intersects $X$ in a curve of genus zero, which must be birational to $\P^1_k$ by Bertini's Theorem because it contains the smooth point $P$.  Let $A$ be the set of linear spaces $L$ that satisfy the following conditions:
\begin{itemize}
\item $L$ is defined over $k$.
\item $L$ contains $P$.
\item $L\cap X$ is a curve of genus zero that is smooth at $P$.
\item $L\cap N$ is empty, so that $\mathcal{L}\cap\mathcal{N}$ is an arithmetic zero-cycle, where $\mathcal{L}$ is the closure of $L$ in $\mathcal{X}$.
\end{itemize}
By Theorem~\ref{curvesintegral} applied to $D'$ and $N$, for any $L\in A$, there are infinitely many rational points $Q$ on $L\cap X$ such that
\[\lambda_{D,v}(Q)\leq n_v\]
for all $v\not\in S$.  (Since $D'$ has degree at most two, the curve $L\cap X$ is a curve of genus zero whose normalization is either $\P^1$, $\mathbb{G}_a$, or $\mathbb{G}_m$.  The first two kinds of curves cannot have a finite but nonempty set of integral points for any $k$, and the infinitude of the unit group of $\O_k$ ensures the same for $\mathbb{G}_m$.)  The union of such curves is dense in $X$.  To see this, note that the union of all linear spaces $L$ that pass through $P$ and are defined over $k$ is Zariski dense in $X$.  The conditions of being smooth at $P$ and intersecting $N$ trivially are both nonempty open conditions (because $N$ has codimension at least two in $X$), so the union of the curves in $A$ must be dense in $X$.  We conclude that there is a Zariski dense set of $(D,S)$-integral points, as desired.  \qed

\end{document}